\documentstyle[12pt]{article}
\begin{document}
\begin{center}
 {\LARGE{}A Riemann--Lebesgue Lemma for Jacobi expansions}
 \\[.4cm]
 {\sc{}George  Gasper\footnote{Department of Mathematics, Northwestern 
 University, Evanston, IL 60208, USA. The work of this author was 
 supported in part by the National Science   Foundation under grant 
 DMS-9401452.}
 and Walter  Trebels\footnote{Fachbereich Mathematik, TH Darmstadt, 
 Schlo\ss{}gartenstr.7, D--64289
 Darmstadt, Germany.} }\\[.4cm]
 {\it Dedicated to P.~L. Butzer on the occasion of his  $65$-th birthday} \\ 
 [.4cm]
 {(Sept. 26, 1994 version)}
 \end{center}

 \bigskip
 {\bf Abstract.}  A Lemma of Riemann--Lebesgue type for Fourier--Jacobi 
coefficients is derived. Via integral representations of Dirichlet--Mehler 
type for Jacobi polynomials its proof directly 
reduces to the classical Riemann--Lebesgue Lemma for Fourier coefficients. 
Other proofs are sketched. Analogous results are also derived for 
Laguerre expansions and for Jacobi transforms.


\bigskip
{\bf AMS(MOS) subject classifications.} 33C45,  42A16, 42B10, 42C10, 44A20

\bigskip

\bigskip \noindent 
The classical Riemann--Lebesgue Lemma states (see \cite[p. 168]{bune}, 
\cite[(4.4), p. 45]{zy}):\\

\smallskip \noindent
{\it If $f\in L^1(-\pi ,\pi ),$ then 
$$\lim _{|k|\to \infty } \int _{-\pi }^\pi f(\theta )e^{-ik\theta } d\theta =0.$$
Also, by using the identities 
$e^{i \theta }e^{i \phi} =  e^{i (\theta + \phi)}$ and 
$e^{i \theta } = \cos \theta + i \sin \theta,$  if 
$a, b, c_0,c_1 \in {\bf R}$ and $f\in L^1(a,b),$  then
\begin{equation}\label{rl-lemma}
\lim _{|k|\to \infty } \int _a^b f(\theta ) \cos \, ((k +c_0)\theta +c_1)\,
d\theta  =0 .
\end{equation} }

\smallskip \noindent
Here we will first consider the extension of this result to 
Fourier--Jacobi coefficients. For this purpose we  need to introduce the 
following notation.
Fix $\alpha, \beta >-1$ and let $L_{(\alpha 
,\beta )}$ denote the space of measurable functions on $[0,\pi ]$ with finite 
norm 
$$\| f \| _{L_{(\alpha, \beta )}} =\int _0^\pi \Big| f(\theta )\Big|  
\Big( \sin \frac{\theta }{2} \Big) ^{2\alpha +1} 
\Big( \cos \frac{\theta }{2} \Big) ^{2\beta +1}d\theta  .$$
Define normalized Jacobi polynomials by 
$R_k^{(\alpha, \beta )} (x)=P_k^{(\alpha, \beta )} (x)/P_k^{(\alpha ,
\beta )} (1)$, where $P_k^{(\alpha, \beta )} (x)$ is the Jacobi polynomial of 
degree $k$ and order $(\alpha, \beta )$, see \cite{szego}. For  $f\in 
L_{(\alpha, \beta )} $, its $k$-th
Fourier--Jacobi coefficient ${\hat f}_{(\alpha, \beta )} (k)$ is defined by
\begin{equation}\label{coeff}
{\hat f}_{(\alpha, \beta )} (k)=\int _0^\pi f(\theta )\,R_k^{(\alpha, \beta )} 
(\cos \theta )\Big( \sin \frac{\theta }{2} \Big) ^{2\alpha +1} 
\Big( \cos \frac{\theta }{2} \Big) ^{2\beta +1} d\theta \, .
\end{equation}
Then $f$ has an expansion of the form
$$f(\theta ) \sim \sum _{k=0}^\infty {\hat f}_{(\alpha, \beta )} (k)\,
h_k^{(\alpha, \beta )}R_k^{(\alpha, \beta )} (\cos \theta )\, ,$$
where the normalizing factors $h_k^{(\alpha, \beta )} $ are given by 
$$h_k^{(\alpha, \beta )} = 
\Big( \|(R_k^{(\alpha, \beta )})^2\| _{L_{(\alpha, \beta )}}\Big) ^{-1} 
\quad \quad \quad \quad \quad \quad \quad \qquad \quad \qquad \qquad \qquad $$
$$\quad \quad \quad \quad \quad =
\frac{(2k+\alpha +\beta +1) \Gamma (k+\alpha +\beta 
+1) \Gamma (k+\alpha +1)}{\Gamma (k+\beta +1) \Gamma (k+1) \Gamma (\alpha +1)
\Gamma (\alpha +1)} \approx (k+1)^{2\alpha +1} $$
and the $\approx $ sign means that there are 
positive constants $C,C'$ such that \break $C'h_k^{(\alpha ,\beta )} \le (k+1)
^{2\alpha +1} \le C h_k^{(\alpha ,\beta )}$ holds.

\smallskip \noindent
In a recent paper \cite{revisited} on ultraspherical multipiers we used the
$\alpha=\beta$ case of the observation that if $(\alpha, \beta) \in S$ with
$$ S:= \{(\alpha, \beta): \alpha \ge \beta >-1, \alpha \ge -1/2 \}, $$ 
then, by \cite[Theorem 7.32.1]{szego}, 
$$\max_{-1\le x \le 1} |R^{(\alpha, \beta )}_k(x)| \le 1,\qquad k\in {\bf N_0},$$ 
and hence
$$|{\hat f}_{(\alpha, \beta )} (k)| \le \| f \| _{L_{(\alpha, \beta )}},
\qquad k\in {\bf N_0}.$$

\smallskip \noindent
Since 
$R^{(-1/2 ,-1/2 )}_k(\cos \theta) = \cos k\theta$ and thus, 
 by  (\ref{rl-lemma}), 
$\lim _{k\to \infty  } {\hat f}_{(-1/2 ,-1/2)} (k)=0$
for $f\in L_{(-1/2 ,-1/2)}, $ this led us to consider
the problem of determining all $(\alpha, \beta)$ with  $\alpha, \beta >-1$
for which
\begin{equation}\label{rljacobi} 
\lim _{k\to \infty  } {\hat f}_{(\alpha, \beta )} (k)=0
\end{equation}
\smallskip \noindent
whenever $f\in L_{(\alpha, \beta )}, $ and to derive the following extension
of the  Riemann--Lebesgue Lemma  to  Fourier--Jacobi coefficients.

\medskip

\bigskip \noindent
{\bf Lemma.} 
{\it Let $\alpha, \beta >-1$. Then (\ref{rljacobi}) holds for each 
$f\in L_{(\alpha, \beta )} $ if and only if $(\alpha, \beta) \in S.$}

\medskip 

\bigskip \noindent
{\bf Proof.} 
First consider the case $\beta > \alpha >-1$. Formula 16.4 (1) in 
\cite[p. 284]{tit} (corrected by inserting a missing $n!$ factor into the 
denominator of the right hand side) shows for
$\alpha > -1, \; \beta+\rho > -1$ that
\begin{equation}\label{counterex}
|\int _{-1}^1 R^{(\alpha, \beta )}_k(x) (1-x)^\alpha 
(1+x)^{\beta+\rho} dx|  \approx (k+1)^{-2\rho-\alpha -\beta-2} \, .
\end{equation}
 Since $-2\rho-\alpha -\beta-2 = \beta - \alpha -2(\beta+\rho+1),$
it follows that if $\beta > \alpha > -1$ and $0<\beta+\rho+1<(\beta-\alpha)/2,$
then the right hand side of (\ref{counterex}) tends to $\infty$ as $k \to 
\infty.$ 

\medskip \noindent
Now let $-1 < \beta \le  \alpha < -1/2$. Introduce the linear functional 
$T_k :L_{(\alpha , \beta )} \to {\bf C},\;  T_kf:= {\hat f}_{(\alpha, \beta )}
(k).$
By \cite[Theorem 7.32.1]{szego},  $|P^{(\alpha, \beta )}_k(x)|$ attains its 
maximum at a point $x'$ ( one of the two maximum points
nearest $x_0=(\beta -\alpha )/(\alpha +\beta +1)$) and 
$|P^{(\alpha, \beta )}_k(x')| \sim (k+1)^{-1/2}$. By the continuity of the 
Jacobi polynomial $P^{(\alpha, \beta )}_k(x)$ there exists a $\delta _k >0$ 
such that   
$2|P^{(\alpha, \beta )}_k(x)| \ge |P^{(\alpha, \beta )}_k(x')|$ for all $x$ 
with $|x-x'| \le \delta _k  \, .$
Now choose $f_k \in L_{(\alpha , \beta )}$ with 
${\rm supp} f_k \subset \{ x :|x-x'| \le \delta _k \},$ 
$ {\rm sgn} \, f_k(x) ={\rm sgn} \, P^{(\alpha, \beta )}_k(x) \; {\rm if} \; 
|x-x'| \le \delta _k  ,\; \| f_k \| _{L_{(\alpha , \beta )}} =1.$
Then 
$$|T_kf|= |{\hat f}_{\alpha, \beta }(k)| \ge \frac{C|P^{(\alpha, \beta )}_k(x')
|}{2|P^{(\alpha, \beta )}_k(1)|} \,\| f_k \| _{L_{(\alpha , \beta )}} 
\ge C' (k+1)^{-\alpha -1/2} ,$$
i.e., $\| T_k\| \ge C'(k+1)^{-\alpha -1/2} $.
Hence, by the uniform boundedness principle there exists an $ 
f^* \in L_{(\alpha , \beta )}$ with $\lim _{k\to \infty } |T_k f^*| 
= \infty $.  Summarizing,  (\ref{rljacobi}) cannot hold for each 
$f\in L_{(\alpha, \beta )} $
when $(\alpha, \,  \beta) \not\in S.$

\medskip \noindent
Now let $\alpha \ge -1/2,\; \alpha >\beta > -1 $. By the definition of the 
Fourier--Jacobi coefficients, 
$${\hat f}_{(\alpha, \beta )} (k)=\bigg( \int _0^{\pi /2} + \int _{\pi /2}^\pi 
\bigg) f(\theta )R_k^{(\alpha, \beta )} 
(\cos \theta )\Big( \sin \frac{\theta }{2} \Big) ^{2\alpha +1} 
\Big( \cos \frac{\theta }{2} \Big) ^{2\beta +1} d\theta =:I_k +J_k.$$
Then $J_k$ tends to $0$ for $k\to \infty $ if $\alpha >-1/2$ and $\alpha >\beta 
>-1$ since, by \cite[ paragraph below Theorem 7.32.1]{szego}, 
$$\max_{\pi/2 \le \theta \le \pi}|R_k^{(\alpha , \beta)}(\cos \theta)| 
\le C (k+1)^{{\rm max}\{ \beta ,-1/2\} -\alpha}.$$
If $\alpha = \beta$, then $J_k$ is of the same type as $I_k$. Thus we can 
restrict ourselves to a discussion of $I_k$.

\bigskip \noindent
To estimate $I_k$ we will use the formula of Dirichlet--Mehler 
type in Gasper \cite[(6)]{mehler}
\begin{eqnarray}\label{dirmehler}
R_k^{(\alpha, \beta )} (\cos \theta ) & = & \frac{2^{(\alpha +\beta +1)/2}
\Gamma(\alpha+1)}{\Gamma(1/2)\Gamma(\alpha+1/2)} \,
 (1-\cos \theta )^{-\alpha } \nonumber\\
 {} & {} &  \times \,  
\int_0^\theta \cos \,(k+(\alpha +\beta +1)/2)\phi \, \frac{(\cos \phi-\cos \theta 
)^{\alpha -1/2}}{(1+\cos \phi )^{(\alpha +\beta )/2}} \nonumber\\
{} & {} & \; \times \,  {}_2F_1\Big[ \frac{\alpha +\beta +1}{2}, 
\frac{\alpha +\beta }{2}; 
\alpha + \frac{1}{2}; \frac{\cos \phi -\cos \theta }{1+\cos \phi } \Big] d\phi,
\end{eqnarray}
which is valid for $\alpha > -1/2, \; 0 < \theta < \pi.$
Inserting this integral representation into $I_k$ and interchanging
 the order of integration we find that 
$$ I_k = C \int _0^{\pi /2} g(\phi) \cos \, (k+(\alpha +\beta +1)/2)\phi \, d\phi$$
with 
$$ g(\phi) =  \Big( \cos \frac{\phi}{2} \Big) ^{-\alpha-\beta }
 \int _\phi ^{\pi /2} f(\theta)  
(\cos \phi - \cos \theta ) ^{\alpha -1/2} 
\Big( \sin \frac{\theta}{2}\Big)
\Big( \cos \frac{\theta}{2} \Big) ^{2\beta +1}$$
$$\quad \quad \quad \quad \quad \quad \quad \quad \quad 
\times \, {}_2F_1\Big[ \frac{\alpha +\beta +1}{2}, \frac{\alpha +\beta }{2}; 
\alpha + \frac{1}{2}; \frac{\cos \phi -\cos \theta }{1+\cos \phi }
 \Big] \,d\theta .$$
Now notice that $0<\cos \theta <\cos \phi<1$ and
$(\cos \phi -\cos \theta )/(1+\cos \phi) < 1/2$ when $0<\phi<\theta <\pi /2.$
Thus the ${}_2F_1$ function in the above integrand is uniformly bounded and
$$\int _0^{\pi /2} |g(\phi)| \, d\phi \le C  \int _0^{\pi /2} 
\int _\phi ^{\pi /2}|f(\theta)|
(\cos \phi - \cos \theta ) ^{\alpha -1/2} 
\Big( \sin \frac{\theta}{2}\Big)
\Big( \cos \frac{\theta}{2} \Big) ^{2\beta +1}\, d\theta \, d\phi$$
$$\le C \int_0^{\pi /2}|f(\theta)|\, h(\theta)
\Big( \sin \frac{\theta }{2} \Big) ^{2\alpha +1} 
\Big( \cos \frac{\theta }{2} \Big) ^{2\beta +1}\, d\theta \qquad $$
with
$$h(\theta) =(\sin \theta )^{-2\alpha } 
\int _0^\theta (\cos \phi -\cos \theta )^{\alpha -1/2} \, d\phi. $$
Since
\begin{eqnarray*}
h(\theta) &  \le  &
C (\sin \theta )^{-2\alpha} \int _0^\theta  \Big( \sin \frac{\theta 
+\phi }{2} \, \sin \frac{\theta -\phi }{2} \Big) ^{\alpha -1/2}\, d\phi \\
{} &  \le & C (\sin \theta )^{-\alpha -1/2} \bigg[ \int _0^{\theta /2} + 
\int _{\theta /2} ^\theta \bigg] \Big( \sin \frac{\theta -\phi }{2} \Big) 
^{\alpha -1/2}\, d\phi \\
{} & \le & C (\sin \theta )^{-1} \int _0^{\theta /2} d\phi \,
+ \,C (\sin \theta )^{ -\alpha -1/2} \int _{\theta /2} ^\theta (\theta -\phi 
)^{\alpha -1/2} d\phi \le C  
\end{eqnarray*}
when $\alpha >-1/2$ and $0< \theta \le \pi /2$, it follows that 
$g(\phi)$ is integrable on $[0, \pi/2]$ and hence $I_k \rightarrow 0$
as $k \rightarrow \infty$ by (\ref{rl-lemma}).

\medskip \noindent
It remains to consider the case
 $\alpha = -1/2, -1 < \beta < -1/2 .$  To handle this case we first
observe that by letting $\alpha$ decrease to $-1/2$ in  (\ref{dirmehler})
and proceeding as in the derivation of formula (6.13) in
\cite{banach} we obtain that
\begin{eqnarray}\label{limdirmehler}
R_k^{(-1/2, \beta )} (\cos \theta ) & = & 
\Big( \cos \frac{\theta}{2} \Big)^{-\beta -1/2} \cos \,(k+\beta/2 +1/4) \theta
\nonumber\\
{} & {} & + \frac{1}{4} (\beta^2-\frac{1}{4}) \Big( \sin \frac{\theta}{2} \Big)
\int_0^\theta \Big( \cos \frac{\phi}{2} \Big)^{-\beta -3/2} 
\cos \,(k+\beta/2 +1/4)\phi  \nonumber\\
{} & {} & \; \; \; \; \times \, {}_2F_1\Big[ \beta/2 +5/4,\beta/2 +3/4; 2; 
\frac{\cos \phi -\cos \theta }{1+\cos \phi } \Big] \,d\phi 
\end{eqnarray}
for $0 < \theta < \pi .$ Since the series 
${}_2F_1[ \beta/2 +5/4,\beta/2 +3/4; 2; x]$
 converges at $x=1$ when $\beta < 0,$ it converges uniformly on $[0, 1],$ 
and thus,  observing that 
$$ 0< \frac{\cos \phi -\cos \theta }{1+\cos \phi } <1\, ,\quad \quad 
0< \phi <\theta < \pi ,$$ 
it is clear that the 
${}_2F_1$ series in the above integral is uniformly bounded when $\beta < 0.$ 
Hence, for $-1 < \beta < -1/2$, the use of (\ref{limdirmehler}) in 
(\ref{coeff}) gives 
$$ {\hat f} _{(-1/2,\beta )}(k) = c_1 \int _0^\pi f(\theta ) \Big( \cos 
\frac{\theta }{2} \Big) ^{\beta +1/2} \cos \,(k+\beta /2+1/4)\theta \, d\theta 
\quad \quad \quad \quad \quad \quad \quad $$
$$+c_2 \int _0^\pi f(\theta ) \Big( \sin \frac{\theta}{2} \Big) \Big( \cos 
\frac{\theta }{2} \Big) ^{2\beta +1} \bigg[ \int _0^\theta  \Big( \cos 
\frac{\phi }{2}\Big) ^{-\beta -3/2} \cos \,(k+\beta/2 +1/4)\phi $$
$$ \quad \quad \quad \quad \times \, {}_2F_1\Big[ \beta/2 +5/4,\beta/2 +3/4; 2; 
\frac{\cos \phi -\cos \theta }{1+\cos \phi } \Big] d\phi \, \bigg] \, d\theta 
= M_k+N_k ,$$
say. Since  $\| f\| _{L_{(-1/2,\beta /2 -1/4)}} 
\le \| f\| _{L_{(-1/2,\beta )}}$ for $\beta \le -1/2,$ $M_k \rightarrow 0$
as $k \rightarrow \infty$ by the Riemann--Lebesgue Lemma (\ref{rl-lemma}).

\smallskip \noindent
Concerning $N_k$, after an interchange of integration  one arrives at 
$$N_k = c_2 \int _0^\pi \cos \,(k+\beta/2 +1/4)\phi \,\bigg[ \Big( \cos 
\frac{\phi }{2}\Big) ^{-\beta -3/2} \int _\phi ^\pi
f(\theta ) \Big( \sin \frac{\theta}{2} \Big) \Big( \cos \frac{\theta }{2} \Big) 
^{2\beta +1} $$
$$ \quad \quad \quad \quad \times \, {}_2F_1\Big[ \beta/2 +5/4,\beta/2 +3/4; 2; 
\frac{\cos \phi -\cos \theta }{1+\cos \phi } \Big] d\theta \, \bigg] \, d\phi 
$$
and the assertion again follows by (\ref{rl-lemma}) once 
we have shown that $\int _0^\pi |[\ldots ]|\, d\phi < \infty .$ But this is 
immediate, since the occurring ${}_2F_1$ function is uniformly bounded,  
$$ \Big| \int _\phi ^\pi \ldots \, d\theta \Big| \le C \int _0^\pi |f(\theta )| 
\Big( \cos \frac{\theta }{2} \Big) ^{2\beta +1} d\theta < \infty,
\qquad \qquad 0 < \phi < \pi, $$
 and 
$$    \int _0^\pi \Big( \cos \frac{\phi }{2}\Big) ^{-\beta -3/2} d\phi < \infty 
\, , \quad\quad \beta < -1/2.$$
Thus the Lemma is established.

\bigskip \noindent
{\bf Remarks.} 1)  Since 
$P_k^{(\alpha, \beta )} (-x) = (-1)^k P_k^{(\beta, \alpha)} (x)$ and
$P_k^{(\alpha, \beta )} (1) = $ $k+\alpha \choose k$ $ \approx (k+1)^\alpha,$
the  Lemma implies that if $f\in L_{(\alpha, \beta )} $ 
and $\max \{\alpha, \beta \} \ge -1/2,$ then
$$\int _0^\pi f(\theta )P_k^{(\alpha, \beta )} 
(\cos \theta )\Big( \sin \frac{\theta }{2} \Big) ^{2\alpha +1} 
\Big( \cos \frac{\theta }{2} \Big) ^{2\beta +1} d\theta \, 
= o(k^{\max \{\alpha, \beta \}}), \qquad k \rightarrow \infty.$$

\medskip \noindent
2) Notice that the above proof via the Dirichlet-Mehler type integral is an 
elementary one; it reduces the problem straight to the classical 
Riemann-Lebesgue Lemma for Fourier coefficients and does not use 
any density properties of subspaces of $L_{(\alpha , \beta )}$. By making use 
of such properties we can give the following simpler proofs. It is well known 
that the subspaces of cosine polynomials, simple functions, and of step 
functions are dense in $L_{(\alpha , \beta 
)}.$ Thus, if $f\in L_{(\alpha, \beta )} $  and 
 $\varepsilon > 0,$ then we can write $f=g+h,$ where 
$\| h\| _{L_{(\alpha ,\beta )}} < \varepsilon$ and 
$g$ is a cosine polynomial, a simple function, or a step function. 
Now let $(\alpha ,\beta )\in S.$ Since
$$ |{\hat h}_{(\alpha , \beta )}(k)|\le 
\| h\| _{L_{(\alpha ,\beta )}} < \varepsilon ,\qquad  k\in {\bf N_0},$$
because $|R_k^{(\alpha ,\beta )}(\cos \, \theta)| \le 1$ for  
$(\alpha ,\beta )\in S,$ to prove that 
${\hat f}_{(\alpha , \beta )}(k) \to 0 $ as $k\to 
\infty $ it suffices to show that ${\hat g}_{(\alpha , \beta )}(k) \to 0 $ 
as $k\to \infty. $ 
This is obvious when $g$ is a cosine
polynomial since ${\hat g}_{(\alpha , \beta )}(k)=0$ for all sufficiently large 
$k.$ 
 If $g$ is a simple function then, being bounded, 
 it is square integrable with respect to the weight function
$\Big( \sin \frac{\theta }{2} \Big) ^{2\alpha +1} 
\Big( \cos \frac{\theta }{2} \Big) ^{2\beta +1}$ and the Parseval formula gives 
$\sum (k+1)^{2 \alpha +1} |{\hat g}_{(\alpha , \beta )}(k)|^2 < \infty $, which
implies that ${\hat g}_{(\alpha , \beta )}(k) \to 0 $ as $k\to 
\infty $ when $(\alpha ,\beta )\in S.$ If $g$ is a step function then
it is a finite linear combination of characteristic functions $\chi(\theta)$
of subintervals of $ (0, \pi),$ so that it remains to show for such
$\chi$ that
${\hat \chi}_{(\alpha , \beta )} (k) \to 0 $ as $k\to 
\infty;$ but this easily follows by using the integral 
\cite[10.8 (38)]{htf} and the asymptotic expansion \cite[(8.21.10)]{szego}.

\medskip \noindent
In the case of Laguerre expansions one does not have a Dirichlet--Mehler type 
formula at one's disposal. However, the preceding three arguments apply. 
To sketch this 
we introduce the  Lebesgue space 
$$L_{w(\alpha )} = \{ f: \; \| f\| _{L_{w(\alpha )}} =  \int _0^\infty
|f(x)| e^{-x/2} x^{\alpha }\, dx   < \infty \} \; ,  \quad \quad
 \alpha >-1, $$
and the normalized Laguerre polynomials $R_k^\alpha (x)$ 
$$ R_k^\alpha (x)=L_k^\alpha (x)/L_k^\alpha (0),\quad \; \; \; 
 L_k^\alpha (0)=A_k^\alpha =
\left( \matrix{ k+\alpha \cr k} \right) =\frac{\Gamma (k+ \alpha +1)}{\Gamma
(k+1) \Gamma (\alpha +1)  }\, , $$
where  $L_k^\alpha (x),  k\in {\bf N}_0$,
is the classical Laguerre polynomial of degree $k$ and
order $\alpha$ (see Szeg\"o \cite[p. 100]{szego}).
Associate to $f$ its formal Laguerre series
$$ f(x) \sim (\Gamma (\alpha +1))^{-1}
\sum _{k=0}^\infty {\hat f}_\alpha (k) L_k^\alpha (x), $$
where the Fourier Laguerre coefficients of $f$ are defined by
$${\hat f}_\alpha (k) = \int _0^\infty f(x) R_k^\alpha (x) x^\alpha e^{-x}\, dx 
$$
when the integrals exist ( for a more detailed description see, 
e.g., \cite{laguerre}). If $\alpha \ge 0$ then, by \cite[10.18 (14)]{htf},
$$| e^{-x/2} R_k^{\alpha}(x)| \le 1\, , \quad \quad x \ge 0, \quad
 k\in {\bf N_0}.$$
Since polynomials and simple functions are dense in $L_{w(\alpha )},$ by
proceeding as above it follows again that 
 a Riemann-Lebesgue Lemma for Fourier-Laguerre coefficients 
holds:
$$ \lim _{k\to \infty } {\hat f}_\alpha (k) =0,\quad f\in L_{w(\alpha )} \, , 
\quad \alpha \ge 0\, .$$
This can also be proved by using the
density of step functions in $L_{w(\alpha )}$
and the observation that, by \cite[10.12 (28)]{htf} and \cite[(8.22.1)]{szego},
$$\int_0^a R_k^{\alpha} (x) e^{-x} x^\alpha \, dx =
\frac{1}{\alpha+1} \, e^{-a} a^{\alpha+1} R_{k-1}^{\alpha+1}(a) \to 0,  \qquad
k\to \infty, $$ 
for each $a>0.$ 

\medskip \noindent
3)  One could ask: Does a Riemann--Lebesgue Lemma also hold for 
the system  
\break
$ \{ \sqrt{h_k^{(\alpha ,\beta )}}
R_k^{(\alpha ,\beta )}(\cos \theta) \},$ 
which is orthonormal with respect to the weight function
\break
$( \sin \frac{\theta }{2} ) ^{2\alpha +1} 
( \cos \frac{\theta }{2}) ^{2\beta +1}?$ 
That this cannot be true for $\alpha > -1/2$
  can be seen by an argument analogous to that
at the beginning of the proof of the Lemma: introduce a corresponding linear 
functional, estimate its norm from below by considering a neighborhood of 
$x=1$, and apply the uniform boundedness principle. This also shows (replace 
$(k+1)^{\alpha +1/2}$ by some $(k+1)^\varepsilon ,\; \varepsilon >0 )$ that 
a ``better'' result than that given in the Lemma, better in the sense 
that for general $f\in L_{(\alpha ,\beta )}$ there is a certain rate of 
decrease of the Fourier--Jacobi coefficients ${\hat f}_{(\alpha ,\beta )}(k)$, 
cannot hold.

\medskip \noindent
4) For Fourier coefficients of a function $f\in L^1(-\pi,\pi) $ 
it is well known that they decrease faster for smoother functions 
(see, e.g., \cite[(4.3), p. 45]{zy}). This phenomenon also occurs for Jacobi 
expansions. Let us illustrate this by considering a special case of Besov 
spaces  investigated by Runst and Sickel \cite{rusi} (for  
$\alpha \ge \beta \ge -1/2 $): 
Let $ \delta >0$. We say that $f\in B^\delta _{1,\infty ,\alpha, \beta }$ if 
$$  \| f\| _{B^\delta _{1,\infty ,\alpha, \beta }} :=
\sup _{j\in {\bf N}_0} 2^{j\delta } \Big\| \sum _{k=0}^\infty 
{\hat f}_{(\alpha, \beta )} (k) \,\varphi _j(k) \,
h_k^{(\alpha, \beta )} R_k^{(\alpha, \beta )} (\cos \theta ) \Big\| 
_{L_{(\alpha, \beta )}} < 
\infty \, ,$$
where $\varphi (x) \in C^\infty ({\bf R})$ has compact support in $(2^{j-1}, 
2^{j+2})$ and is identicallly $1$ for $2^j \le x \le 2^{j+1}$. 
Then for $n\in  [2^j ,2^{j+1}]$ one obtains 
$$ |{\hat f}_{(\alpha, \beta )} (n) | \le 
\Big\| \sum _{k=0}^\infty {\hat f}_{(\alpha, \beta )} (k) \, \varphi _j(k) \,
h_k^{(\alpha, \beta )} R_k^{(\alpha, \beta )} (\cos \theta ) \Big\| 
_{L_{(\alpha, \beta )}} \le C n^{-\delta } 
\| f\| _{B^\delta _{1,\infty ,\alpha, \beta }} $$
uniformly in $j$.
Bavinck \cite{ba} introduced Lipschitz spaces based on the generalized Jacobi 
translation operator (see \cite{banach}). These coincide with the above Besov
spaces (see \cite[p. 374]{ba} and \cite[Remark 15 and Theorem 5]{rusi} 
and observe that the domain 
of the infinitesimal generator considered by Bavinck is just the domain of the 
square of the infinitesimal generator considered by Runst and Sickel). 

\smallskip \noindent
5) In the same spirit we can extend the Riemann--Lebesgue Lemma for 
cosine transforms on the  half-axis to Jacobi transforms. For $\alpha >-1$ and 
$\beta \in {\bf R}$ we denote the underlying 
space of measurable functions on ${\bf R}_+$ by
$$L_{(\alpha ,\beta )}({\bf R}_+) = \{ f\, :\; \| f\| _{L_{(\alpha ,\beta )}
({\bf R}_+)}:= \int _0^\infty |f(t)| \Big( \sinh t\Big) ^{2\alpha +1}   
\Big( \cosh t\Big) ^{2\beta +1} dt < \infty \} \, .$$ 
Then the Jacobi transform 
of an $L_{(\alpha ,\beta )}({\bf R}_+) $ function $f$ is defined by 
$$ {\cal J}^{(\alpha ,\beta )}[f](\tau ) = \frac{2^{2(\alpha +\beta +1)+1/2}}
{\Gamma (\alpha +1)} \int _0^\infty f(t) \, \varphi _\tau ^{(\alpha ,\beta )}(t) 
\Big( \sinh t\Big) ^{2\alpha +1} \Big( \cosh t\Big) ^{2\beta +1} dt ,$$
whenever the integral converges, where $\varphi _\tau ^{(\alpha ,\beta )}(t)$
is the Jacobi function  defined by
$$ \varphi _\tau ^{(\alpha ,\beta )}(t) = {}_2F_1\Big[ \frac{1}{2}(\rho 
+i\tau ), \frac{1}{2}(\rho -i\tau ); \alpha +1; -(\sinh t)^2 \Big] \, $$
with $\rho =\alpha + \beta +1$; see, e.g., Koornwinder \cite{ko}. 
Note that ${\cal J}^{(-1/2,-1/2)}[f]$ is just 
the cosine transform of $f$, see \cite[(3.4)]{ko}. 
This time, one can reduce the problem to the classical Riemann-Lebesgue Lemma 
for the cosine transform: for $\alpha > -1/2$ 
Koornwinder \cite[(2.21)]{ko} has shown the following 
integral representation of Dirichlet--Mehler type
$$\varphi _\tau ^{(\alpha ,\beta )}(t) = 2^{-\alpha +3/2} 
\frac{\Gamma (\alpha +1)}{\Gamma (\alpha +1/2)\Gamma (1/2)} 
\frac{1}{(\sinh t)^{2\alpha }(\cosh t)^{\alpha +\beta}} \quad \quad \quad \quad 
\quad \quad \quad \quad \quad $$
$$ \times  \int _0^t \cos \tau s \; 
(\cosh \, 2t -\cosh \, 2s )^{\alpha - 
1/2} {}_2F_1\Big[ \alpha +\beta ,\alpha -\beta ;\alpha +\frac{1}{2}; 
\frac{\cosh \, t -\cosh \, s}{2\cosh \, t} \Big] \, ds\, .$$
 From this integral it follows (cf. \cite[p. 150]{ko}) that
\begin{equation}\label{estimate}
|\varphi _\tau ^{(\alpha ,\beta )}(t)| \le C (1+t) e^{-(\alpha +\beta +1 )t}
\, ,\quad t, \tau \in {\bf R}_+ \, , \quad \alpha > -1/2.
\end{equation}
Hence the Jacobi transform of a function
$f\in L_{(\alpha ,\beta )}({\bf R}_+) $ exists 
as a uniformly bounded function of $\tau \in {\bf R}_+ $ if
$\alpha > -1/2$ and $\alpha + \beta > -1.$
Then, by proceeding as in proof of the Lemma, 
the  Riemann-Lebesgue Lemma for cosine transforms now implies that 
$ {\cal J}^{(\alpha ,\beta )}[f](\tau )$ vanishes at infinity when
$\alpha > -1/2$ and $\alpha + \beta > -1.$  This result 
can also be proved by using
the density in $L_{(\alpha, \beta)} ({\bf R_+})$ of finite linear
combinations of characteristic functions of bounded intervals, 
\, \cite[(2.10)]{ko}, \, (\ref{estimate}), and the method described in Remark 2.
This result and the inequality in (\ref{estimate})
can be extended to 
$\alpha = -1/2, \,\alpha +\beta > -1$  by using the $\alpha \searrow -1/2$ limit case
of the above integral representation:
$$\varphi _\tau ^{(-1/2 ,\beta )}(t) = 
(\cosh \, t)^{-\beta -1/2} \cos \,\tau t \qquad \qquad \qquad \qquad \qquad \qquad 
\qquad \qquad \qquad $$
$$\qquad \qquad \quad
+ \,(\frac{1}{4}-\beta^2) (\sinh \, t) ( \cosh \,t)^{-\beta -1/2}
 \int _0^t \cos \,\tau s \; 
(\cosh \, t +\cosh \, s )^{-1}$$
$$ \times  \, {}_2F_1\Big[\frac{1}{2} + \beta, \frac{1}{2} - \beta; 2; 
\frac{\cosh \, t -\cosh \, s}{2\cosh \, t} \Big] \, ds\,.$$

\noindent
We intend to consider the general complex parameter case and related 
problems in another paper.

\end{document}